\documentclass[12pt]{article}
\usepackage{amssymb}
\usepackage{amsmath}
\usepackage{graphicx}
\begin{document}

\newtheorem{theorem}{Theorem}

\title{Unfolding a Codimension-Two, Discontinuous, Andronov-Hopf Bifurcation}
\author{D.J.W.~Simpson and J.D.~Meiss\thanks{
D.~J.~W.~Simpson and J.~D.~Meiss gratefully
acknowledge support from NSF grant DMS-0707659.
}\\
Department of Applied Mathematics\\
University of Colorado\\
Boulder, CO 80309-0526}
\maketitle

\begin{abstract}
We present an unfolding of the codimension-two scenario
of the simultaneous occurrence of a discontinuous bifurcation
and an Andronov-Hopf bifurcation in a
piecewise-smooth, continuous system of
autonomous ordinary differential equations in the plane.
We find the Hopf cycle undergoes a grazing bifurcation
that may be very shortly followed by a saddle-node bifurcation
of the orbit.
We derive scaling laws for
the bifurcation curves that emanate from the
codimension-two bifurcation.
\end{abstract}

\section{Introduction}

A system of differential equations is said to be {\em piecewise-smooth}
if it is everywhere smooth except on some codimension-one boundaries
called {\em switching manifolds}.
This paper is concerned with such systems that are continuous
everywhere but non-differentiable on switching manifolds.
Piecewise-smooth continuous systems have been utilized to model
a range of diverse physical situations,
for instance vibro-impacting mechanical systems \cite{LeNi04,WiDe00},
switching in electrical circuits \cite{BaVe01,ZhMo03,Ts03}
and various non-smooth phenomena in biology and physiology \cite{Ro70,KeSn98}.

The interaction of invariant sets with switching manifolds
often produces bifurcations not seen in smooth systems.
The last two decades have seen an explosion of interest in such bifurcations
and many new results, see for instance
\cite{DiBu08,LeNi04,ZhMo03} and references within.
In the neighborhood of a single switching manifold a piecewise-smooth, continuous system
may be written as
\begin{equation}
\dot{x} = \left\{ \begin{array}{lc}
F^{(L)}(x), & H(x) \le 0 \\
F^{(R)}(x), & H(x) \ge 0
\end{array} \right.\;,
\label{eq:genSys}
\end{equation}
where $x \in \mathbb{R}^N$,
$F^{(L)},F^{(R)} : \mathbb{R}^N \to \mathbb{R}^N$ are $C^k$ and
$H : \mathbb{R}^N \to \mathbb{R}$ is sufficiently smooth.
The switching manifold is the set
$\mathcal{S} = \{ x \in \mathbb{R}^N ~|~ H(x) = 0 \}$
and by continuity, $F^{(L)} = F^{(R)}$ on $\mathcal{S}$.

A point $x^* \in \mathbb{R}^N$ is an equilibrium of
the left-half-system if $F^{(L)}(x^*) = 0$ and is said
to be {\em admissible} if $H(x^*) \le 0$ and {\em virtual} otherwise
(and vice-versa for the right-half-system).
By continuity of (\ref{eq:genSys}), an equilibrium that lies exactly on $\mathcal{S}$
is an admissible equilibrium of both smooth half-systems.
This codimension-one phenomenon generally gives rise
to what is known as a {\em discontinuous bifurcation}
(sometimes called a {\em boundary equilibrium bifurcation})
for which there are two basic generic scenarios.
Either equilibria in the two half-systems coexist, collide and annihilate
on the switching manifold ({\em non-smooth saddle-node})
or one equilibrium is admissible on each side
of the bifurcation ({\em persistence}).
In addition, other invariant sets may be created at the bifurcation,
such as a periodic orbit in a manner akin to an Andronov-Hopf
bifurcation \cite{FrPo97,SiMe07}.

Even though linear terms of an appropriate series expansion of the system
generally completely determine dynamical behavior
local to a discontinuous bifurcation,
in a general $N$-dimensional system the bifurcation may be extremely complex.
Much work has been done for low ($N \le 3$) dimensional systems
\cite{FrPo98,CaFr05,LeNi04,DiBu08}.
In two-dimensions, if a periodic orbit emanates from a generic discontinuous
bifurcation, it must encircle an equilibrium of focus type \cite{FrPo98}.
In this paper we unfold about the codimension-two point where the focus-type
equilibrium has purely imaginary eigenvalues at the crossing point.
This scenario corresponds to the simultaneous occurrence of a discontinuous
bifurcation and a (smooth) Andronov-Hopf bifurcation
and has recently been observed in a model of yeast growth \cite{SiKo08}.

We find that as parameters are changed to move away from the Hopf bifurcation,
the associated Hopf cycle grows in size as usual,
until grazing the switching manifold.
No bifurcation occurs at the grazing point in the sense
that the phase portrait does not change topologically
(because the system is continuous at the switching manifold \cite{DiBu01}).
However, very shortly beyond the grazing bifurcation
a saddle-node bifurcation of the orbit may occur.
We derive a condition governing when this occurs
and a scaling law describing how close to the grazing
the saddle-node bifurcation occurs.
Our results are presented in Theorem \ref{th:c22dhb}.

The remainder of this paper is organized as follows.
In \S\ref{sec:prelim} we transform the general system
to a normal form
involving companion matrices and state the theorem.
In \S\ref{sec:example} we provide a simple example
and use it to illustrate the theorem
and numerically verify predicted scaling laws.
\S\ref{sec:outline} outlines our method of proof
and \S\ref{sec:proof} presents a proof of the theorem.
Conclusions are presented in \S\ref{sec:conclusions}.

\section{Preliminaries and Theorem Statement}
\label{sec:prelim}

Consider a two-dimensional, piecewise-$C^k$ continuous system
of ordinary differential equations in $\mathbb{R}^2$
with two independent parameters, $\mu$ and $\eta$.
For our analysis below we will need to assume $k \ge 8$.
In a neighborhood of a single switching manifold the system may be written as
\begin{equation}
\left[ \begin{array}{c}
\dot{x} \\ \dot{y}
\end{array} \right] =
\left\{ \begin{array}{lc}
\left[ \begin{array}{c}
f^{(L)}(x,y;\mu,\eta) \\ g^{(L)}(x,y;\mu,\eta)
\end{array} \right], & H(x,y;\mu,\eta) \le 0 \\
\left[ \begin{array}{c}
f^{(R)}(x,y;\mu,\eta) \\ g^{(R)}(x,y;\mu,\eta)
\end{array} \right], & H(x,y;\mu,\eta) \ge 0
\end{array} \right.\;,
\label{eq:genSystem}
\end{equation}
where $H:\mathbb{R}^4 \to \mathbb{R}$ is a sufficiently smooth
(at least $C^3$) function.
The switching manifold is the parameter dependent set,
$\mathcal{S}_{\mu,\eta} = \{ (x,y)^{\sf T} ~|~ H(x,y;\mu,\eta) = 0 \}$.
Without loss of generality we may assume that $H$ vanishes
at $(x,y;\mu,\eta) = (0,0;0,0)$.
If in addition
$\frac{\partial H}{\partial x} + \frac{\partial H}{\partial y} \ne 0$,
then locally $\mathcal{S}_{0,0}$ is a curve intersecting the origin.
Via coordinate transformations in a similar manner as
given in \cite{DiBu01}, we may assume to order $C^3$ that
$H$ is simply equal to $x$.
The higher order terms in $H$ do not affect our analysis below;
thus in what follows, we will assume
$H$ is identically equal to $x$.
The switching manifold is then simply the $y$-axis
and we will refer to
$(f^{(L)},g^{(L)})^{\sf T}$ as the {\em left-half-system} and 
$(f^{(R)},g^{(R)})^{\sf T}$ as the {\em right-half-system}.

We may assume
that there is a discontinuous bifurcation at the origin
when $\mu = \eta = 0$.
Since the origin lies on the switching manifold
and (\ref{eq:genSystem}) is continuous, it is an
equilibrium of both the left and right-half-systems.
In this paper we are interested in the scenario that the
equilibrium in the left-half-plane has
complex-valued eigenvalues
$\lambda_\pm = \nu \pm {\rm i} \omega$.
(We make no assumptions about eigenvalues of the
equilibrium solution in the right-half-plane.)
Assume that when $\mu = \eta = 0$, the eigenvalues are purely imaginary, i.e.,
\begin{equation}
\nu(0,0) = 0,~~\omega(0,0) > 0\;.
\end{equation}
Notice $\pm {\rm i} \omega(0,0)$ are the eigenvalues of the matrix 
\begin{equation}
J = \left[ \begin{array}{cc}
\frac{\partial f^{(L)}}{\partial x} & \frac{\partial f^{(L)}}{\partial y} \\
\frac{\partial g^{(L)}}{\partial x} & \frac{\partial g^{(L)}}{\partial y}
\end{array} \right] \Bigg|_{(x,y;\mu,\eta) = (0,0;0,0)}\;.
\end{equation}
Therefore, in particular
\begin{eqnarray}
\det(J) & = & \omega^2(0,0) \ne 0 \;, \label{eq:detNonZero} \\
J_{12} & = & \frac{\partial f^{(L)}}{\partial y} \bigg|_{(0,0;0,0)}
\ne 0\;. \label{eq:topLeftNonZero}
\end{eqnarray}
By the implicit function theorem and (\ref{eq:detNonZero})
the left-half-system has an equilibrium
$(x^{*(L)}(\mu,\eta),y^{*(L)}(\mu,\eta))^{\sf T}$
where $x^{*(L)}$ and $y^{*(L)}$ are $C^k$ functions and
$x^{*(L)}(0,0) = y^{*(L)}(0,0) = 0$.
As is generically the case,
we assume the distance of the equilibrium from the switching manifold
varies linearly with some combination of parameters.
Without loss of generality we may assume $\mu$ is a suitable choice.
That is
\begin{equation}
\frac{\partial x^{*(L)}}{\partial \mu}(0,0) \ne 0\;.
\label{eq:xiDotBNonZero}
\end{equation}
Again by the implicit function theorem, there is a $C^k$ function,
$\phi$, such that $x^{*(L)}(\phi(\eta),\eta) = 0$.
In other words when $\mu = \phi(\eta)$,
the equilibrium lies on the switching manifold.
After performing the nonlinear change of coordinates
\begin{equation}
\begin{split}
\mu & \mapsto \mu - \phi(\eta)\;, \\
y & \mapsto y - y^{*(L)}(\phi(\eta),\eta)\;,
\label{eq:fTrans1}
\end{split}
\end{equation}
we may factor $\mu$ out of the constant term in the system (\ref{eq:genSystem}),
i.e.\footnote{
We use $O(k)$ (and $o(k)$)
to denote terms that are order $k$ (larger than order $k$)
in all variables and parameters.
When necessary to distinguish orders
we are more specific, e.g. $O(|x,y|^3)$.
}
\begin{equation}
\left[ \begin{array}{c}
f^{(L)}(0,0;\mu,\eta) \\ g^{(L)}(0,0;\mu,\eta)
\end{array} \right] =
\left[ \begin{array}{c}
f^{(R)}(0,0;\mu,\eta) \\ g^{(R)}(0,0;\mu,\eta)
\end{array} \right] =
\left[ \begin{array}{c}
p(\mu,\eta) \\ q(\mu,\eta)
\end{array} \right] \mu + o(k)\;,
\end{equation}
where $p$ and $q$ are $C^{k-1}$.
The left and right-half-systems may now be written as
\begin{eqnarray}
\left[ \begin{array}{c}
f^{(i)}(x,y;\mu,\eta) \\ g^{(i)}(x,y;\mu,\eta)
\end{array} \right] & = &
\left[ \begin{array}{c}
p(\mu,\eta) \\ q(\mu,\eta)
\end{array} \right] \mu +
\left[ \begin{array}{cc}
a_i(\mu,\eta) & b(\mu,\eta) \\
c_i(\mu,\eta) & d(\mu,\eta)
\end{array} \right]
\left[ \begin{array}{c}
x \\ y
\end{array} \right] \nonumber \\
& & +~\left[ \begin{array}{c}
\tilde{f}^{(i)}(x,y;\mu,\eta) \\ \tilde{g}^{(i)}(x,y;\mu,\eta)
\end{array} \right] + o(k)\;,
\label{eq:semiRedForm}
\end{eqnarray}
where $a_i$, $b$, $c_i$, $d$ are each $C^{k-1}$
and $\tilde{f}^{(i)}$,$\tilde{g}^{(i)}$ are $C^k$ functions
nonlinear in $x$ and $y$.
(As with $p$ and $q$, subscripts are not required for the coefficients
$b$ and $d$ because the system is continuous.)

We show now that it is possible to change coordinates to set $p(\mu,\eta) \equiv 0$.
Since $b(0,0) \ne 0$ by (\ref{eq:topLeftNonZero}),
using the implicit function theorem,
there exists a unique $C^k$ function, $\psi$,
such that $f^{(L)}(0,\psi(\mu,\eta);\mu,\eta) = 0$.
After the coordinate change
\begin{equation}
y \mapsto y - \psi(\mu,\eta)\;,
\label{eq:pToZero}
\end{equation}
the system remains in the form (\ref{eq:semiRedForm})
with $p(\mu,\eta) \equiv 0$.

Finally we transform the system to the usual canonical companion matrix form.
As is well-known,
this may be accomplished when the system
is {\em observable} in the control theory sense \cite{CaFr02,DiBu08}.
Our system is observable by (\ref{eq:topLeftNonZero}),
and the required transformation is
\begin{equation}
\begin{split}
y & \mapsto - d(\mu,\eta) x + b(\mu,\eta) y\;, \\
\mu & \mapsto -b(\mu,\eta) q(\mu,\eta) \mu\;, \\
\eta & \mapsto a_L(\mu,\eta) + d(\mu,\eta)\;.
\label{eq:hatCoord}
\end{split}
\end{equation}
The first equation is nonsingular by (\ref{eq:topLeftNonZero}).
The second equation is nonsingular by (\ref{eq:xiDotBNonZero})
and the final equation is nonsingular if we make a final nondegeneracy assumption:
\begin{equation}
\frac{\partial \nu}{\partial \eta}(0,0) \ne 0\;.
\end{equation}
The transformed system has a normal form suited for further analysis
and is given in the following theorem.

\begin{figure}[h]
\begin{center}
\includegraphics[width=15cm,height=6cm]{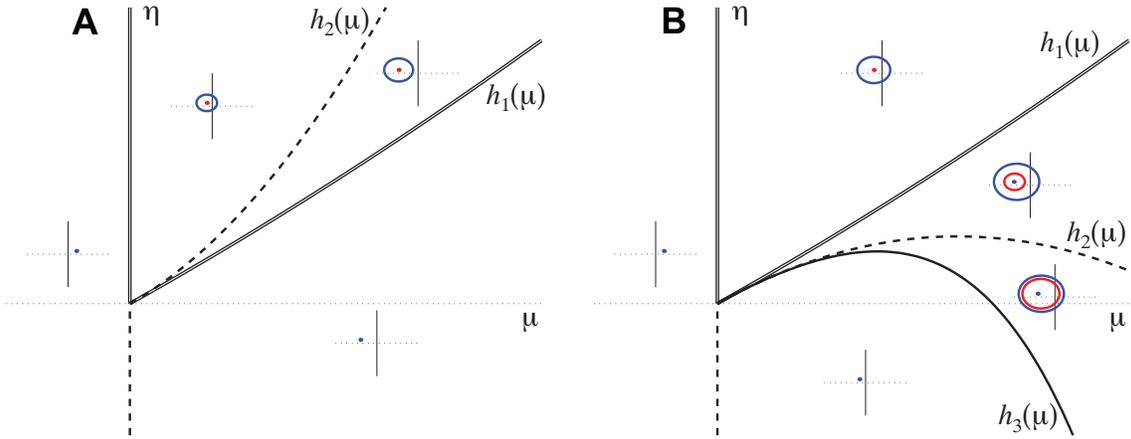}
\caption{Schematics showing unfoldings predicted by Theorem \ref{th:c22dhb}
when $\tau_R < 0$, $\delta_R > 0$.
In panel A, $a_0 < 0$, and in panel B, $a_0 > 0$.
Double lines correspond to Hopf bifurcations or discontinuous bifurcations
that create a periodic orbit.
Dashed lines correspond to a collision of an equilibrium or periodic
orbit with the switching manifold at which no bifurcation occurs.
Solid lines correspond to saddle-node bifurcations of periodic orbits.
Included are phase portraits showing local behavior.
Stable orbits are colored blue;
unstable orbits are colored red.
Thin solid lines denote the switching manifold.
\label{fig:unfoldingSchHb}}
\end{center}
\end{figure}

\begin{theorem}~\\
Consider the two-dimensional, piecewise-$C^k$ ($k \ge 8$),
continuous system of differential equations
\begin{equation}
\left[ \begin{array}{c}
\dot{x} \\ \dot{y}
\end{array} \right] =
\left\{ \begin{array}{lc}
\left[ \begin{array}{c}
f^{(L)}(x,y;\mu,\eta) \\ g^{(L)}(x,y;\mu,\eta)
\end{array} \right], & x \le 0 \\
\left[ \begin{array}{c}
f^{(R)}(x,y;\mu,\eta) \\ g^{(R)}(x,y;\mu,\eta)
\end{array} \right], & x \ge 0
\end{array} \right.\;,
\label{eq:theSystem}
\end{equation}
where
\begin{equation}
\left[ \begin{array}{c}
f^{(i)}(x,y;\mu,\eta) \\ g^{(i)}(x,y;\mu,\eta)
\end{array} \right] =
\left[ \begin{array}{c}
0 \\ -\mu
\end{array} \right] +
\left[ \begin{array}{cc}
\tau^{(i)}(\mu,\eta) & 1 \\
-\delta^{(i)}(\mu,\eta) & 0
\end{array} \right]
\left[ \begin{array}{c}
x \\ y
\end{array} \right] +
\left[ \begin{array}{c}
\tilde{f}^{(i)}(x,y;\mu,\eta) \\ \tilde{g}^{(i)}(x,y;\mu,\eta)
\end{array} \right] + o(k)\;,
\label{eq:canonFormHB}
\end{equation}
such that
\begin{equation}
\tau^{(L)}(\mu,\eta) \equiv \eta\;,
\end{equation}
and $\tilde{f}^{(i)}$,$\tilde{g}^{(i)}$ consist of terms
that are nonlinear in $x$ and $y$.
Suppose,
\begin{enumerate}
\renewcommand{\labelenumi}{\roman{enumi})}
\item $\delta^{(L)}(0,0) = \omega^2$ for $\omega > 0$,
\item $a_0 \ne 0$, where
\end{enumerate}
\vspace{-5mm}
\begin{eqnarray}
a_0 & = & \frac{1}{16}(f^{(L)}_{xxx} + g^{(L)}_{xxy} +
\omega^2 f^{(L)}_{xyy} + \omega^2 g^{(L)}_{yyy}) -
\frac{1}{16} f^{(L)}_{xy}(f^{(L)}_{xx} + \omega^2 f^{(L)}_{yy}) \nonumber \\
& & +~\frac{1}{16} g^{(L)}_{xy}(\frac{1}{\omega^2} g^{(L)}_{xx} + g^{(L)}_{yy}) +
\frac{1}{16}(\frac{1}{\omega^2}f^{(L)}_{xx}g^{(L)}_{xx} -
\omega^2 f^{(L)}_{yy}g^{(L)}_{yy})
\label{eq:a0HB}
\end{eqnarray}
evaluated at $(x,y;\mu,\eta) = (0,0;0,0)$, and
\begin{enumerate}
\renewcommand{\labelenumi}{\roman{enumi})}
\item[iii)] $\tau_R = \tau^{(R)}(0,0) \ne 0$.
\end{enumerate}
Then, near $(x,y;\mu,\eta) = (0,0;0,0)$,
if $\delta_R = \delta^{(R)}(0,0) \ne 0$
when $\delta_R \mu < 0$
there is a unique equilibrium in
the right-half-plane given by $C^k$ functions
\begin{equation}
\left[ \begin{array}{c}
x^{*(R)}(\mu,\eta) \\ y^{*(R)}(\mu,\eta)
\end{array} \right] =
\left[ \begin{array}{c}
-\frac{1}{\delta_R} \\ \frac{\tau_R}{\delta_R}
\end{array} \right] \mu + O(2)\;.
\end{equation}
The equilibrium is repelling if $\tau_R, \delta_R > 0$,
attracting if $-\tau_R, \delta_R > 0$
and a saddle if $\delta_R < 0$.
When $\mu > 0$, there is a unique equilibrium
in the left-half-plane given by $C^k$ functions
\begin{equation}
\left[ \begin{array}{c}
x^{*(L)}(\mu,\eta) \\ y^{*(L)}(\mu,\eta)
\end{array} \right] =
\left[ \begin{array}{c}
-\frac{1}{\omega^2} \\ 0
\end{array} \right] \mu + O(2)\;.
\end{equation}
Furthermore, there exist unique $C^{k-1}$, $C^{k-1}$, $C^{k-2}$ functions
$h_1,h_2,h_3 : \mathbb{R} \to \mathbb{R}$ respectively, with
\begin{eqnarray}
h_1(\mu) & = &
\frac{1}{\omega^2}(f^{(L)}_{xx}+g^{(L)}_{xy})\big|_{(0,0;0,0)} \mu + O(\mu^2)\;,
\label{eq:h1HB} \\
h_2(\mu) & = &
h_1(\mu) - \frac{2 a_0}{\omega^4} \mu^2 + O(\mu^3)\;,
\label{eq:h2HB} \\
h_3(\mu) & = &
h_2(\mu) - \frac{8 \pi^2 a_0^3}{3 \omega^{12} \tau_R^2} \mu^6 + o(\mu^6)\;,
\label{eq:h3HB}
\end{eqnarray}
such that when $\mu > 0$,
\begin{enumerate}
\renewcommand{\labelenumi}{\roman{enumi})}
\item the curve $\eta = h_1(\mu)$ corresponds to a locus of Andronov-Hopf bifurcations
of $(x^{*(L)},y^{*(L)})^{\sf T}$ that are supercritical if $a_0 < 0$
and subcritical if $a_0 > 0$; this equilibrium
is attracting if $\eta < h_1(\mu)$ and repelling if $\eta > h_1(\mu)$,
\item the curve $\eta = h_2(\mu)$ corresponds to a locus of grazing bifurcations
of the associated Hopf cycle with the $y$-axis,
\item if $a_0 \tau_R < 0$, the Hopf cycle exists for values of $\eta$
between $h_1(\mu)$ and $h_3(\mu)$, a periodic orbit of opposing stability exists for
$\eta < h_3(\mu)$ if $a_0 < 0$ and $\eta > h_3(\mu)$ if $a_0 > 0$
and the two orbits coincide at a locus of saddle-node bifurcations
$\eta = h_3(\mu)$,
\item if $a_0 \tau_R > 0$, the Hopf cycle exists for
$\eta > h_1(\mu)$ if $a_0 < 0$ and $\eta < h_1(\mu)$ if $a_0 > 0$.
\end{enumerate}
\label{th:c22dhb}
\end{theorem}
Theorem \ref{th:c22dhb} predicts essentially two different unfolding scenarios.
These are illustrated in Fig.~\ref{fig:unfoldingSchHb}.

\section{Example}
\label{sec:example}
As an example consider the piecewise-$C^\infty$, continuous system
\begin{equation}
\begin{split}
\dot{u} & = -\alpha + \frac{2}{15} \beta + v + \frac{1}{5} u^2 + u^3\;, \\
\dot{v} & = -\frac{5}{4} \alpha + \frac{1}{6} \beta - \frac{3}{8} u +
\frac{1}{10}(\beta - 1) v + |\frac{1}{8} u - \frac{1}{10} v|\;.
\label{eq:exampleHB}
\end{split}
\end{equation}
When $\alpha = \beta = 0$, the origin is an equilibrium on a
switching manifold.
Its two one-sided limiting associated eigenvalues are
$\pm \frac{1}{\sqrt{2}} {\rm i}$ and
$-\frac{1}{10} \pm \frac{\sqrt{6}}{5} {\rm i}$.
This example exhibits the codimension-two scenario in which we are interested
and satisfies the required non-degeneracy conditions.

\begin{figure}[h]
\begin{center}
\includegraphics[width=10.8cm,height=9cm]{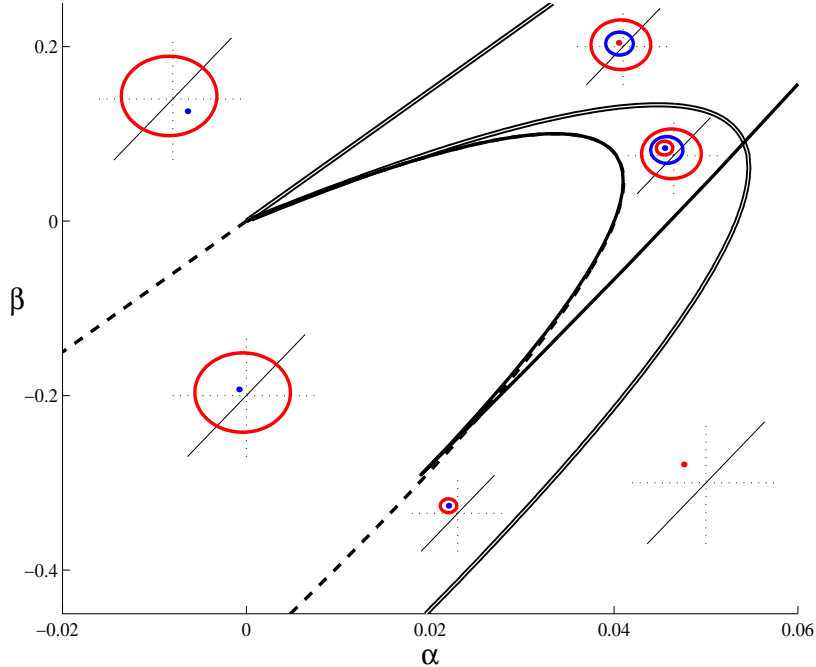}
\caption{
A bifurcation set of the example system (\ref{eq:exampleHB}).
Sketches of representative phase portraits are included.
The line styles and color scheme is the same as in Fig.~\ref{fig:unfoldingSchHb}.
\label{fig:exBifSet}}
\end{center}
\end{figure}

A bifurcation set for (\ref{eq:exampleHB}) is shown in Fig.~\ref{fig:exBifSet}.
For small values of $\alpha$ and $\beta$
the bifurcation set is a smooth distortion of
Fig.~\ref{fig:unfoldingSchHb}, panel B.
However here there also exists an unstable periodic orbit.
The predictions of Theorem \ref{th:c22dhb} break down away from
$\alpha = \beta = 0$ where this orbit
collides with the stable orbit.
The resulting locus of saddle-node bifurcations of periodic orbits intersects the
saddle-node locus anticipated by Theorem \ref{th:c22dhb}
at a cusp bifurcation at $(\alpha,\beta) \approx (0.019,-0.29)$.

In order to compare the predictions of Theorem \ref{th:c22dhb} and
in particular the scaling laws (\ref{eq:h1HB})-(\ref{eq:h3HB}),
we transform the system (\ref{eq:exampleHB})
to the form given in the theorem.
The switching manifold is the line $u = \frac{4}{5} v$,
therefore we let
\begin{equation}
x = u - \frac{4}{5} v\;.
\end{equation}
This particular example has been chosen because the
transformation may be computed explicitly.
We have $\phi(\beta) = -\frac{5}{3} \beta$,
$y^{*(L)}(\phi(\beta),\beta) = 0$ and
$\psi(\alpha,\beta) = 0$.
Combining the individual transformations (\ref{eq:fTrans1}), (\ref{eq:pToZero}) and
(\ref{eq:hatCoord}) produces
\begin{equation}
\begin{split}
y & = -\frac{1}{10}(\beta - 4) u + v\;, \\
\mu & = -\frac{1}{10} (\frac{33}{2} -
\beta)(\alpha - \frac{2}{15} \beta)\;, \\
\eta & = \frac{1}{10} \beta\;.
\end{split}
\end{equation}
In the transformed coordinates the system is
\begin{equation}
\left[ \begin{array}{c} \dot{x} \\ \dot{y} \end{array} \right] =
\left\{ \begin{array}{lc}
\left[ \begin{array}{c} 0 \\ -\mu \end{array} \right] +
\left[ \begin{array}{cc} \eta & 1 \\ -\frac{1}{2} & 0 \end{array} \right]
\left[ \begin{array}{c} x \\ y \end{array} \right] +
\left[ \begin{array}{c} 1 \\ \frac{2}{5}-\eta \end{array} \right]
(\frac{1}{5} u^2 + u^3), & x \le 0 \\
\left[ \begin{array}{c} 0 \\ -\mu \end{array} \right] +
\left[ \begin{array}{cc} \eta-\frac{1}{5} & 1 \\ -\frac{1}{4} & 0 \end{array} \right]
\left[ \begin{array}{c} x \\ y \end{array} \right] +
\left[ \begin{array}{c} 1 \\ \frac{2}{5}-\eta \end{array} \right]
(\frac{1}{5} u^2 + u^3), & x \ge 0
\end{array} \right.\;,
\label{eq:transExHB}
\end{equation}
where
\begin{equation}
u = \frac{25 x + 20 y}{33 - 20 \eta}\;.
\end{equation}
The values of the important constants are
\begin{equation}
\begin{split}
(a_0,\omega,\tau_R,\delta_R) & =
\left( \frac{25}{88},\frac{1}{\sqrt{2}},-\frac{1}{5},\frac{1}{4} \right)\;, \\
(f^{(L)}_{xx},g^{(L)}_{xy}) \big|_{(0,0;0,0)} & =
\left( \frac{250}{1089},\frac{80}{1089} \right)\;.
\end{split}
\end{equation}

\begin{figure}[h]
\begin{center}
\includegraphics[width=15cm,height=4cm]{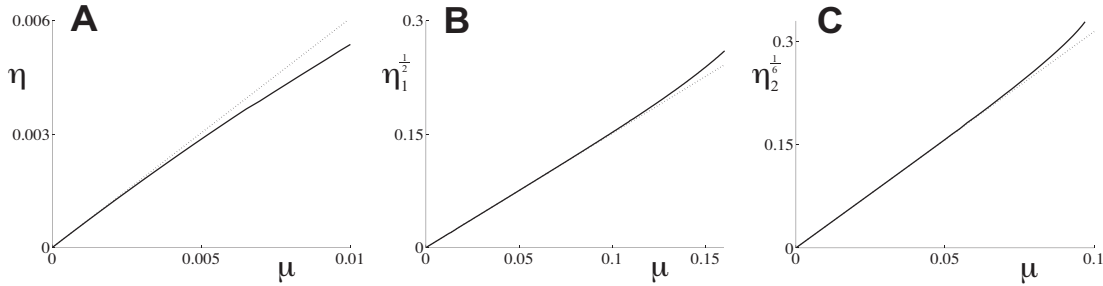}
\caption{Numerical verification of Theorem \ref{th:c22dhb}
for the transformed example system (\ref{eq:transExHB}).
Panel A shows $h_1(\eta)$.
Panel B is a plot of the square root of the difference
between $h_2(\eta)$ and $h_1(\eta)$.
Panel C is a plot of the sixth root of the difference
between $h_3(\eta)$ and $h_2(\eta)$.
The dashed lines in the three panels are the lowest order scaling
predictions (\ref{eq:h1HB})-(\ref{eq:h3HB}), respectively.
\label{fig:scVer}}
\end{center}
\end{figure}

Fig.~\ref{fig:scVer} shows comparisons of numerical computations of
the curves $h_1$, $h_2$ and $h_3$, with their predicted scalings.
We find the numerical results are in full agreement with Theorem \ref{th:c22dhb}.
A calculation of error terms to quantitatively estimate the difference
between the lowest order approximations and the true curves
is beyond the scope of this paper.

\section{Proof Outline}
\label{sec:outline}

\begin{figure}[h]
\begin{center}
\includegraphics[width=10.8cm,height=9cm]{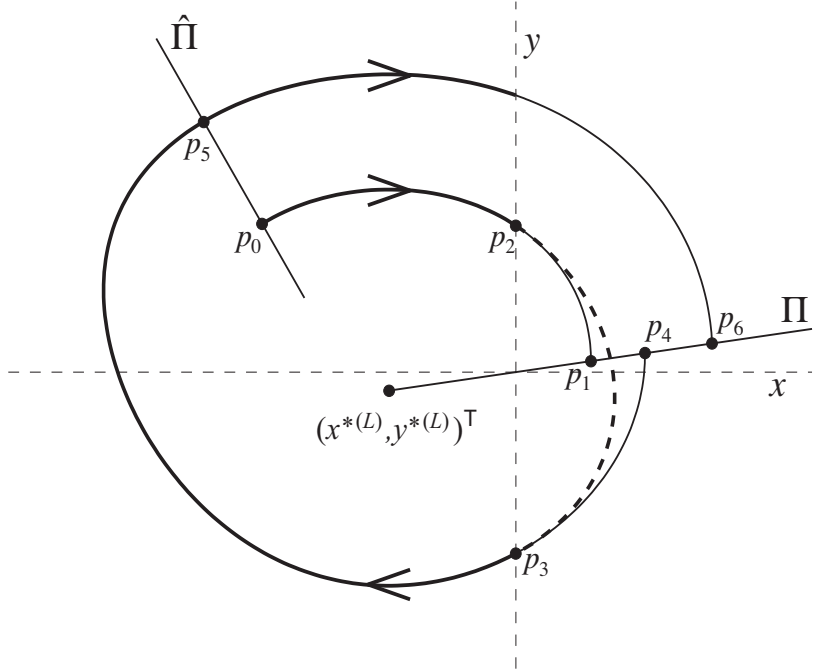}
\caption{
Schematic illustrating a trajectory and Poincar\'{e} sections
relating to Theorem \ref{th:c22dhb} when $\mu > 0$.
A solid [dotted] curve corresponds to following the
left [right]-half-flow.
The true trajectory is the thick curve.
The Poincar\'{e} maps constructed in the proof are
$p_2 = \mathcal{P}_1(p_1)$,
$p_3 = \mathcal{P}_2(p_2)$,
$p_4 = \mathcal{P}_3(p_3)$,
$p_6 = \mathcal{P}_{\rm lhf}(p_4)$.
\label{fig:qz}
}
\end{center}
\end{figure}

Here we present an outline for the proof of Theorem \ref{th:c22dhb}
that follows in \S\ref{sec:proof}.
By assumption, when $\mu > 0$, there is an equilibrium,
$(x^{*(L)},y^{*(L)})^{\sf T}$, in the left-half-plane.
Close to this point there will exist an appropriate Poincar\'{e} section, $\hat{\Pi}$,
in the left-half-plane, see Fig.~\ref{fig:qz}.
The trajectory that begins from a point $p_0$ on $\hat{\Pi}$,
spirals clockwise around $(x^{*(L)},y^{*(L)})^{\sf T}$
and reintersects $\hat{\Pi}$ at some point $p_5$.
We are interested in periodic orbits of the flow,
thus when $p_5 = p_0$.

It will be more convenient to use a second Poincar\'{e} section, $\Pi$, that is a 
semi-infinite line intersecting the equilibrium and
the origin.
Artificially following the left-half-flow from $p_0$,
we arrive at a point $p_1$ on $\Pi$.
If $p_1$ were in the left-half-plane,
by continuing from $p_1$ along the left-half-flow we would arrive at $p_5$.
However if $p_1$ is in the right-half-plane,
we must first calculate a correction, following
the left-half-flow from a different point on $\Pi$, $p_4$,
in order to arrive at $p_5$.
The correction results from the lack of smoothness at the switching manifold.
The mapping between $p_1$ and $p_4$ is known as the discontinuity map
\cite{DaNo00,DiBu01}.

To compute the discontinuity map,
we follow the left-half-flow back from $p_1$ until we arrive at a point, $p_2$,
on the $y$-axis after a time $T_1 < 0$.
We then follow the right-half-flow from $p_2$ until the next intersection
with the $y$-axis at $p_3$ after a time $T_2 > 0$.
Finally we follow the left-half-flow to $p_4$ on $\Pi$ after a time $T_3 < 0$.
Let $\mathcal{P}_1$, $\mathcal{P}_2$ and $\mathcal{P}_3$ denote maps
relating to these three steps, respectively.
The discontinuity map is then
\begin{equation}
\mathcal{P}_{\rm dm} = \mathcal{P}_3 \circ
\mathcal{P}_2 \circ \mathcal{P}_1\;.
\end{equation}
As will be shown, $\mathcal{P}_{\rm dm}$ has a $\frac{3}{2}$-type singularity.
This is because the system is continuous on the switching manifold \cite{DiBu01,DiBu08}.
When $p_1$ lies in the left-half-plane the discontinuity map
is taken to be the identity map.

A map from $p_0$ to $p_5$ may be derived by composing a map
from $\hat{\Pi}$ to $\Pi$ with $\mathcal{P}_{\rm dm}$ and with 
a map from $\Pi$ to $\hat{\Pi}$.
However, a map from $p_1$ to $p_6$ is simpler
(where $p_6$ is the point on $\Pi$ obtained by following the left-half-flow
from $p_5$) and equivalent.
Let $\mathcal{P}_{\rm lhf}$ denote the map following the left-half-flow
from $\Pi$ to itself and let $T_4$ denote the corresponding
transition time.
Then $p_1$ is mapped to $p_6$ via
\begin{equation}
\mathcal{P} = \mathcal{P}_{\rm lhf} \circ \mathcal{P}_{\rm dm}\;.
\end{equation}
A fixed point of $\mathcal{P}$ corresponds to $p_6 = p_1$.
Notice $p_6 = p_1$ exactly when $p_5 = p_0$.

The most intensive computations in the following proof
are in deriving the maps $\mathcal{P}_{\rm lhf}$,
$\mathcal{P}_1$, $\mathcal{P}_2$ and $\mathcal{P}_3$ described above.
Once this is accomplished, $\mathcal{P}$ is obtained by composition
and it remains to derive its fixed points to find $h_2$
and obtain the function $h_3$ by 
locating where $\mathcal{P}$ has a saddle-node bifurcation.
In order to derive expressions,
not for $h_2$ and $h_3$,
but rather for $h_2 - h_1$ and $h_3 - h_2$,
we will introduce adjusted parameters $\eta_1$ and $\eta_2$,
that represent the deviation from $h_1$ and $h_2$ respectively.

\section{Proof of Theorem \ref{th:c22dhb}}
\label{sec:proof}

For ease of notation we expand the vector field in the left-half-plane
to second order in $x$ and $y$
\begin{equation}
\begin{split}
f^{(L)}(x,y;\mu,\eta) & =
\eta x + y + a_1 x^2 + a_2 xy + a_3 y^2 + O(|x,y|^3)\;, \\
g^{(L)}(x,y;\mu,\eta) & =
-\mu - \delta^{(L)}x + b_1 x^2 + b_2 xy + b_3 y^2 + O(|x,y|^3)\;, \label{eq:fLgLexpanded} 
\end{split}
\end{equation}
where coefficients vary with $\mu$ and $\eta$.

\noindent
{\bf Step 1:}~~~{\em Compute equilibria, eigenvalues and $h_1$.}\\
By the implicit function theorem, there exists a unique $C^k$ function,
$\left[ \begin{array}{c} x^{*(L)} \\ y^{*(L)} \end{array} \right]
: \mathbb{R}^2 \to \mathbb{R}^2$, such that
$\left[ \begin{array}{c}
f^{(L)}(x^{*(L)}(\mu,\eta),y^{*(L)}(\mu,\eta);\mu,\eta) \\
g^{(L)}(x^{*(L)}(\mu,\eta),y^{*(L)}(\mu,\eta);\mu,\eta)
\end{array} \right] = 0$, for small $\mu$, $\eta$.
Via substitution of a series expansion of $x^{*(L)}$ and $y^{*(L)}$
into (\ref{eq:fLgLexpanded})
it is readily determined that
\begin{equation}
\begin{split}
x^{*(L)}(\mu,\eta) & =
-\frac{1}{\omega^2} \mu + \left( \frac{b_1}{\omega^6}
+ \frac{\delta_\mu^{(L)}(0,0)}{\omega^4} \right) \mu^2
+ \frac{\delta_\eta^{(L)}(0,0)}{\omega^4} \mu \eta + O(3)\;, \\
y^{*(L)}(\mu,\eta) & =
-\frac{a_1}{\omega^4} \mu^2 + \frac{1}{\omega^2} \mu \eta + O(3)\;, \label{eq:xyStL}
\end{split}
\end{equation}
where $\mu$ and $\eta$ subscripts denote derivatives.
Near $(x,y;\mu,\eta) = (0,0;0,0)$, $(x^{*(L)},y^{*(L)})^{\sf T}$ is the
only equilibrium of the left-half-flow
and by (\ref{eq:xyStL})
it is an equilibrium of the full flow (i.e.~admissible),
exactly when $\mu \ge 0$.
If $\delta_R \ne 0$, the corresponding equilibrium
in the right-half-plane may be determined similarly.

Let $J^{(L)}$ denote the Jacobian of (\ref{eq:theSystem}) when $x < 0$.
For small $\mu$ and $\eta$ the matrix
\begin{equation}
J^{(L)}(x^{*(L)}(\mu,\eta),y^{*(L)}(\mu,\eta);\mu,\eta)
= \left[ \begin{array}{cc}
j_1 & j_2 \\ j_3 & j_4 \end{array} \right]\;,
\label{eq:DFHB}
\end{equation}
has the complex conjugate eigenvalue pair
\begin{equation}
\lambda_\pm = \nu \pm {\rm i} \xi\;,
\nonumber
\end{equation}
where $\nu$ and $\xi$ are $C^{k-1}$ functions of $\mu$ and $\eta$.
To first order we have
\begin{equation}
\nu(\mu,\eta) = -\frac{2 a_1 + b_2}{2 \omega^2} \mu + \frac{1}{2} \eta + O(2)\;.
\label{eq:nuHb}
\end{equation}
By the implicit function theorem,
there exists a unique $C^{k-1}$ function $h_1 : \mathbb{R} \to \mathbb{R}$,
such that $\nu(\mu,h_1(\mu)) = 0$ for small $\mu$.
Furthermore
\begin{equation}
h_1(\mu) = \frac{2 a_1 + b_2}{\omega^2} \mu + O(\mu^2)\;,
\nonumber
\end{equation}
thereby confirming (\ref{eq:h1HB}) of the theorem.
Let
\begin{equation}
\eta_1 = \eta - h_1(\mu)
\label{eq:etaHatDef}
\end{equation}
represent the deviation from the Hopf bifurcation curve, then
\begin{equation}
\nu(\mu,\eta_1) = \eta_1 (\frac{1}{2} + O(1))\;.
\label{eq:nuHb2}
\end{equation}

\noindent
{\bf Step 2:}~~~{\em Compute $a_0$ and introduce polar coordinates.}\\
Let $v_{\lambda_\pm}$ denote the complex-valued eigenvector
associated with $\lambda_\pm$ for (\ref{eq:DFHB}).
In a standard manner, we construct a matrix $P$ using the real and
imaginary parts of $v_{\lambda_\pm}$
\begin{equation}
P = \left[ \begin{array}{cc}
1 & 0 \\ \frac{1}{2 j_2}(-j_1+j_4) &
-\frac{1}{2 j_2} \sqrt{-(j_1-j_4)^2 - 4 j_2 j_3}
\end{array} \right]\;,
\nonumber
\end{equation}
which is well-defined (because $j_2 \ne 0$ for sufficiently
small $\mu,\eta$) and non-singular.
Let
\begin{equation}
\left[ \begin{array}{c} u \\ v \end{array} \right] =
P^{-1} \left( \left[ \begin{array}{c} x \\ y \end{array} \right] -
\left[ \begin{array}{c}
x^{*(L)}(\mu,\eta) \\ y^{*(L)}(\mu,\eta)
\end{array} \right] \right)\;.
\nonumber
\end{equation}
Then the left-half-system in $(u,v)$ coordinates becomes
\begin{equation}
\left[ \begin{array}{c} \dot{u} \\ \dot{v} \end{array} \right] =
\left[ \begin{array}{c}
\check{f}(u,v;\mu,\eta) \\ \check{g}(u,v;\mu,\eta)
\end{array} \right] = D^{(L)} \left[ \begin{array}{c} u \\ v \end{array} \right] + O(2)\;,
\nonumber
\end{equation}
where
\begin{equation}
D^{(L)} = P^{-1} J^{(L)} P = \left[ \begin{array}{cc}
\nu & -\xi \\ \xi & \nu \end{array} \right]\;.
\nonumber
\end{equation}
Letting
\begin{equation}
z = u + {\rm i} v\;,
\label{eq:zDef}
\end{equation}
then
\begin{equation}
\dot{z} = \lambda_+ z + O(2)\;.
\label{eq:zDot}
\end{equation}
Following standard proofs of the Hopf bifurcation 
\cite{GuHo86,Ku04,Gl99},
there exists a two-variable polynomial $\zeta$, comprised of only
quadratic and cubic terms such that the near identity transformation
\begin{equation}
w = z + \zeta(z,\bar{z})\;,
\label{eq:wDefHB}
\end{equation}
removes all quadratic terms and all but one cubic term from the
left-half-system (\ref{eq:zDot})
\begin{equation}
\dot{w} = \lambda_+ w + A w^2 \bar{w} + O(4)\;,
\end{equation}
where $A \in \mathbb{C}$ and
\begin{eqnarray}
\mathcal{R}e(A(\mu,\eta)) & = & \bigg[ \frac{1}{16}
(\check{f}_{uuu} + \check{g}_{uuv} + \check{f}_{uvv} + \check{g}_{vvv}) +
\frac{1}{16 \xi} \left( \check{f}_{uv}(\check{f}_{uu} + \check{f}_{vv}) \right. \nonumber \\
& & \left. -~\check{g}_{uv}(\check{g}_{uu} + \check{g}_{vv}) -
\check{f}_{uu} \check{g}_{uu} + \check{f}_{vv} \check{g}_{vv} \right) \bigg]
\bigg|_{(u,v)=(0,0)}\;.
\label{eq:a0hatHB}
\end{eqnarray}
Notice $\mathcal{R}e(A(\mu,\eta))$ is a $C^{k-3}$ function of $\mu$ and $\eta$.
Let
\begin{equation}
a_0 = \mathcal{R}e(A(0,0))\;.
\nonumber
\end{equation}
Using (\ref{eq:a0hatHB})
and since $P(0,0) = \left[ \begin{array}{cc} 1 & 0 \\ 0 & -\omega \end{array} \right]$,
it follows that $a_0$ appears as given in the statement of the theorem, (\ref{eq:a0HB}).

To prove that generic Hopf bifurcations occur along $\mu = h_1(\eta)$
it remains to verify the non-degeneracy conditions of the
Hopf bifurcation theorem \cite{GuHo86,Ku04,Gl99}.
\begin{enumerate}
\renewcommand{\labelenumi}{\roman{enumi})}
\item By construction,
$J^{(L)}(x^{*(L)}(\mu,h_1(\mu)),y^{*(L)}(\mu,h_1(\mu));\mu,h_1(\mu))$
has purely imaginary eigenvalues,
$\pm {\rm i} \xi = \pm {\rm i} (\omega + O(\mu)) \ne 0$.
\item $\frac{\partial \nu}{\partial \eta}(\mu,h_1(\mu)) =
\frac{1}{2} + O(\mu) \ne 0$.
\item By assumption,
$\mathcal{R}e(A(\mu,h_1(\mu))) = a_0 + O(\mu) \ne 0$.
\end{enumerate}
Therefore if $a_0 < 0$, the curve $\eta = h_1(\mu)$
corresponds to supercritical Hopf bifurcations
and stable periodic orbits exist for small $\eta > h_1(\mu)$.
Conversely, if $a_0 > 0$, the curve $\eta = h_1(\mu)$
corresponds to subcritical Hopf bifurcations
and unstable periodic orbits exist for small $\eta < h_1(\mu)$.
Consequently we have proven (i) of the theorem.

We now introduce polar coordinates.
Let
\begin{equation}
w = r {\rm e}^{{\rm i} \theta}\;.
\label{eq:polarDef}
\end{equation}
In polar coordinates the left-half-system is
\begin{equation}
\begin{split}
\dot{r} & = \nu r + \mathcal{R}e(A) r^3 + O(r^4)\;, \\
\dot{\theta} & = \xi + \mathcal{I}m(A) r^2 + O(r^3)\;.
\end{split}
\end{equation}
Denote the components of the $C^k$ flow by $R(r_0,\theta_0,t;\mu,\eta_1)$
and $\Theta(r_0,\theta_0,t;\mu,\eta_1)$ respectively.
Expressions for $R$ and $\Theta$ may be derived
by expanding each as a series in $r_0$ and computing
coefficients by solving initial value problems.
We obtain
\begin{eqnarray}
R(r_0,\theta_0,t;\mu,\eta_1) & = & {\rm e}^{\nu t} r_0 +
\mathcal{R}e(A) \frac{{\rm e}^{3 \nu t} - {\rm e}^{\nu t}}{2 \nu} r_0^3
+ O(r_0^4)\;, \label{eq:REqHb} \\
\Theta(r_0,\theta_0,t;\mu,\eta_1) & = & \theta_0 + \xi t +
\mathcal{I}m(A) \frac{{\rm e}^{2 \nu t} - 1}{2 \nu} r_0^2 + O(r_0^3)\;. \label{eq:ThetaEqHb}
\end{eqnarray}

\noindent
{\bf Step 3:}~~~{\em Define the Poincar\'{e} section, $\Pi$.}\\
By using (\ref{eq:fLgLexpanded}) it is straightforward to show
\begin{equation}
\lim_{\mu \to 0} \left( \frac{y^{*(L)}}{x^{*(L)}} \right) = -\eta\;.
\nonumber
\end{equation}
Thus the quotient
\begin{equation}
S = \frac{y^{*(L)}}{x^{*(L)}}\;,
\nonumber
\end{equation}
is a well-defined $C^k$ function for small $\mu$, $\eta$.
For our analysis we only need the Taylor series of $S$ to first order,
which from (\ref{eq:xyStL}) is readily found to be
\begin{equation}
S(\mu,\eta) = \frac{a_1}{\omega^2} \mu - \eta + O(2)\;.
\nonumber
\end{equation}
For $\mu > 0$,
$y = S x$ is the line intersecting $(x^{*(L)},y^{*(L)})^{\sf T}$
and the origin. Let
\begin{equation}
\Pi = \left\{ (x,y) ~|~ y = Sx, x \ge x^{*(L)} \right\}\;,
\label{eq:PiDef}
\end{equation}
as in Fig.~\ref{fig:qz}.
In the polar coordinates centered at the equilibrium, $(r,\theta)$,
$\Pi$ is described by a $C^k$ function
\begin{equation}
\theta_\Pi(r;\mu,\eta_1) = \check{\theta} + \rho_1 r + O(r^2)\;,
\label{eq:thetaPiHb}
\end{equation}
where $\check{\theta}$ and $\rho_1$ are coefficients dependent
on $\mu$ and $\eta_1$.
The coefficient $\check{\theta}$ describes the angle of $\Pi$ at the
equilibrium in $(u,v)$ coordinates; $\rho_1$ arises
from the nonlinear coordinate change, (\ref{eq:wDefHB}).
Since explicit forms for these coefficients will not be required,
we do not derive them.\\

\noindent
{\bf Step 4:}~~~{\em Derive $\mathcal{P}_{\rm lhf}$ and compute $h_2$.}\\
We now wish to determine the left-half Poincar\'{e} map, $\mathcal{P}_{\rm lhf}$.
To do this we compute the trajectory of a point,
$(r_0,\theta_\Pi(r_0;\mu,\eta_1))$ on $\Pi$,
and find its next intersection,
$(r_1,\theta_\Pi(r_1;\mu,\eta_1))$.
The transition time is
\begin{equation}
T_4(r_0;\mu,\eta_1) = \frac{2 \pi}{\omega} + O(1) \label{eq:T4hb}\;.
\end{equation}
Substituting (\ref{eq:nuHb2}) and (\ref{eq:T4hb})
into (\ref{eq:REqHb}) produces
\begin{eqnarray}
r_1(r_0;\mu,\eta_1) & = & R(r_0,\theta_\Pi(r_0;\mu,\eta_1),
T_4(r_0;\mu,\eta_1);\mu,\eta_1) \nonumber \\
& = & \left( 1 + \eta_1(\frac{\pi}{\omega} + O(|\mu,\eta_1|^1)) \right) r_0 +
\frac{2 \pi a_0}{\omega} r_0^3 + O(4)\;.
\label{eq:r_1_r_0HB}
\end{eqnarray}
Let $s$ denote the distance from
$(x^{*(L)},y^{*(L)})$ in $(u,v)$, (\ref{eq:zDef}), coordinates.
That is
\begin{equation}
s = |z| = (u^2 + v^2)^\frac{1}{2}\;.
\nonumber
\end{equation}
By (\ref{eq:wDefHB}) and (\ref{eq:polarDef})
\begin{equation}
s(r) = r + \sigma_2 r^2 + O(r^3)\;,
\label{eq:sAsFuncOfrHB}
\end{equation}
for some coefficient $\sigma_2$, determined by
$\zeta$, (\ref{eq:wDefHB}).
By combining (\ref{eq:r_1_r_0HB}) and (\ref{eq:sAsFuncOfrHB})
we are able to obtain $s_1$ as a function of $s_0$,
where $s_0$ is a point on $\Pi$ and $s_1$ is the next point on $\Pi$
\begin{equation}
s_1(s_0;\mu,\eta_1) =
\left( 1 + \eta_1(\frac{\pi}{\omega} + O(|\mu,\eta_1|^1)) \right) s_0 +
\frac{\pi \sigma_2}{\omega} \eta_1 s_0^2 +
\frac{2 \pi a_0}{\omega} s_0^3 + O(4)\;.
\label{eq:sMapHb}
\end{equation}
Let $(\varepsilon,S \varepsilon)^{\sf T}$ be a point on $\Pi$ in $(x,y)$ coordinates.
Let 
\begin{equation}
\check{\varepsilon} = \varepsilon - x^{*(L)}\;.
\label{eq:epsCheckDef}
\end{equation}
Since $\check{\varepsilon}$ is a scalar multiple of $s$
and when $\mu = \eta_1 = 0$ we have $\check{\varepsilon} = s$,
it follows that to third order
the map between $\check{\varepsilon}_0$ and $\check{\varepsilon}_1$
is the same as (\ref{eq:sMapHb}), i.e.
\begin{equation}
\check{\varepsilon}_1(\check{\varepsilon}_0;\mu,\eta_1) =
\left( 1 + \eta_1(\frac{\pi}{\omega} + O(|\mu,\eta_1|^1)) \right)
\check{\varepsilon}_0
+ \frac{\pi \sigma_2}{\omega} \eta_1 \check{\varepsilon}_0^2
+ \frac{2 \pi a_0}{\omega} \check{\varepsilon}_0^3 + O(4)\;.
\label{eq:epsCheckMapHb}
\end{equation}
The system (\ref{eq:theSystem}) has a periodic orbit that
grazes the $y$-axis when
$\check{\varepsilon}_0 = \check{\varepsilon}_1 = -x^{*(L)} \ne 0$.
By substituting this into (\ref{eq:epsCheckMapHb}) and dividing through
by $-x^{*(L)}$ we obtain
\begin{equation}
1 = 1 + \eta_1(\frac{\pi}{\omega} + O(|\mu,\eta_1|^1))
- \frac{\pi \sigma_2}{\omega} \eta_1 x^{*(L)}(\mu,\eta_1)
+ \frac{2 \pi a_0}{\omega} x^{*(L)^{\scriptstyle 2}}(\mu,\eta_1) + O(3)\;.
\label{eq:h2Finder}
\end{equation}
Using (\ref{eq:xyStL}) and the implicit function theorem
we find (\ref{eq:h2Finder}) is satisfied when
$\eta_1 = \hat{h}_2(\mu)$ for a $C^k$ function
\begin{equation}
\hat{h}_2(\mu) = -\frac{2 a_0}{\omega^4} \mu^2 + O(\mu^3)\;.
\nonumber
\end{equation}
Let
\begin{eqnarray}
h_2 & = & h_1 + \hat{h}_2\;, \nonumber \\
\eta_2 & = & \eta_1 - \hat{h}_2(\mu)\;. \label{eq:etaHatHatDef}
\end{eqnarray}
We have therefore derived (\ref{eq:h2HB}) and proven (ii) of the theorem.
We now write $\mathcal{P}_{\rm lhf}$ in terms of
$\varepsilon$, $\mu$ and $\eta_2$.
Combining (\ref{eq:epsCheckDef}), (\ref{eq:epsCheckMapHb})
and (\ref{eq:etaHatHatDef}) produces
\begin{eqnarray}
\mathcal{P}_{\rm lhf}(\varepsilon_0;\mu,\eta_2) & = &
\mu \eta_2 (\frac{\pi}{\omega^3} + O(|\mu,\eta_2|^1)) +
\Big( 1 + \frac{\pi}{\omega} \eta_2 + q_1 \mu^2 + q_2 \mu \eta_2 \nonumber \\
& & +~q_3 \eta_2^2 + O(|\mu,\eta_2|^3) \Big) \varepsilon_0 +
O(\varepsilon_0^2)\;,
\label{eq:PlhfHB}
\end{eqnarray}
where
\begin{equation}
q_1 = \frac{4 \pi a_0}{\omega^5}\;,
\label{eq:q1HB}
\end{equation}
and explicit forms for $q_2$ and $q_3$ will not be required subsequently.\\

\noindent
{\bf Step 5:}~~~{\em Derive the discontinuity map, $\mathcal{P}_{\rm dm}$.}\\
In order to avoid singularities, here we assume $\mu > 0$ and
introduce the spatial scaling
\begin{equation}
\left[ \begin{array}{c} \hat{x} \\ \hat{y} \end{array} \right]
= \frac{1}{\mu} 
\left[ \begin{array}{c} x \\ y \end{array} \right]\;.
\label{eq:spatScaling}
\end{equation}
In the left-half-plane
\begin{equation}
\left[ \begin{array}{c}
\dot{\hat{x}} \\ \dot{\hat{y}}
\end{array} \right] =
\left[ \begin{array}{c}
\hat{f}^{(L)}(\hat{x},\hat{y};\mu,\eta) \\ \hat{g}^{(L)}(\hat{x},\hat{y};\mu,\eta)
\end{array} \right] =
\frac{1}{\mu} \left[ \begin{array}{c}
f^{(L)}(\mu \hat{x},\mu \hat{y};\mu,\eta) \\ g^{(L)}(\mu \hat{x},\mu \hat{y};\mu,\eta)
\end{array} \right]\;.
\label{eq:xyHatDot}
\end{equation}
Let
\begin{eqnarray}
X^{(L)}(\hat{y}_0,t;\mu,\eta) & = & A_1 t + A_2 t^2 + A_3 t^3 + O(t^4)\;,
\label{eq:XLdefHB} \\
Y^{(L)}(\hat{y}_0,t;\mu,\eta) & = & \hat{y}_0 + B_1 t + B_2 t^2 + B_3 t^3 + O(t^4)\;,
\end{eqnarray}
denote the components of the $C^k$ left-half-flow
for an initial condition $(0,\hat{y}_0)$, on the switching manifold.
The parameter dependent coefficients are obtained
by solving (\ref{eq:xyHatDot}) using (\ref{eq:fLgLexpanded})
\begin{equation}
\begin{split}
A_1 & = \hat{y}_0 + a_3 \mu \hat{y}_0^2 + O(\hat{y}_0^3)\;, \\
A_2 & = -\frac{1}{2} + (\frac{1}{2} \eta - a_3 \mu) \hat{y}_0 + O(\hat{y}_0^2)\;, \\
A_3 & = \frac{1}{3} a_3 \mu - \frac{1}{6} \eta + O(\hat{y}_0)\;, \\
B_1 & = -1 + b_3 \mu \hat{y}_0^2 + O(\hat{y}_0^3)\;, \\
B_2 & = (-\frac{1}{2} \delta^{(L)} - b_3 \mu) \hat{y}_0 + O(\hat{y}_0^2)\;, \\
B_3 & = \frac{1}{3} b_3 \mu + \frac{1}{6} \delta^{(L)} + O(\hat{y}_0)\;.
\end{split}
\end{equation}
We now derive $\mathcal{P}_1$, $\mathcal{P}_2$ and $\mathcal{P}_3$
in scaled coordinates, (\ref{eq:spatScaling}), beginning with $\mathcal{P}_3$.
The point $p_4 = \mathcal{P}_3(p_3)$ (see Fig.~\ref{fig:qz}) lies on $\Pi$,
therefore to compute the corresponding transition time, $T_3$,
we solve $Y^{(L)} = S X^{(L)}$ for $t$.
The function
\begin{eqnarray}
G_1(\hat{y}_0,t;\mu,\eta) & = &
Y^{(L)}(\hat{y}_0,t;\mu,\eta) - S X^{(L)}(\hat{y}_0,t;\mu,\eta) \nonumber \\
& = & \hat{y}_0 - t - S \hat{y}_0 t + \frac{1}{2} S t^2 + O(3)\;, \nonumber
\end{eqnarray}
is $C^k$ and by the implicit function theorem
there exists a unique $C^k$ function $T_3 : \mathbb{R}^3 \to \mathbb{R}$ such that
$G_1(\hat{y}_0,T_3(\hat{y}_0;\mu,\eta);\mu,\eta) = 0$.
Furthermore
\begin{equation}
T_3(\hat{y}_0;\mu,\eta) = \hat{y}_0 - \frac{1}{2} S \hat{y}_0^2 + O(\hat{y}_0^3)\;.
\label{eq:T3HB}
\end{equation}
Combining (\ref{eq:XLdefHB}) and (\ref{eq:T3HB}) yields
\begin{equation}
\mathcal{P}_3(\hat{y}_0;\mu,\eta) =
\frac{1}{2} \hat{y}_0^2 
+ \frac{1}{3} (\eta + a_3 \mu) \hat{y}_0^3 + O(\hat{y}_0^4)\;.
\label{eq:P3HB}
\end{equation}
Notice (see Fig.~\ref{fig:qz}) that
$\mathcal{P}_1^{-1}$ is the same as $\mathcal{P}_3$
except the $y$-component of the initial condition has opposite sign.
That is, $\hat{y}_0 = \mathcal{P}_1(\hat{\varepsilon}_0;\mu,\eta)$,
whenever $\hat{\varepsilon}_0 = \mathcal{P}_3(-\hat{y}_0;\mu,\eta)$.
By inverting (\ref{eq:P3HB}) we obtain
\begin{equation}
\mathcal{P}_1(\hat{\varepsilon}_0;\mu,\eta) =
\sqrt{2} \hat{\varepsilon}_0^{\frac{1}{2}}
- \frac{2}{3} (\eta + a_3 \mu) \hat{\varepsilon}_0
+ O(\hat{\varepsilon}_0^{\frac{3}{2}})\;.
\label{eq:P1HB}
\end{equation}
In a similar manner as for $\mathcal{P}_3$ we are able to use a series expansion
of the right-half-flow to determine $\mathcal{P}_2$. We obtain
\begin{equation}
\mathcal{P}_2(\hat{y}_0;\mu,\eta) =
-\hat{y}_0 - \frac{2}{3} (\tau^{(R)} + a_3 \mu) \hat{y}_0^2
+ O(\hat{y}_0^3)\;.
\label{eq:P2HB}
\end{equation}
The discontinuity map is
$\mathcal{P}_{\rm dm} = \mathcal{P}_3 \circ \mathcal{P}_2 \circ \mathcal{P}_1$.
Composition of (\ref{eq:P3HB})-(\ref{eq:P2HB}) produces
\begin{equation}
\hat{\varepsilon}_4 =
\mathcal{P}_{\rm dm}(\hat{\varepsilon}_1;\mu,\eta) =
\hat{\varepsilon}_1 + \frac{4 \sqrt{2}}{3}
(\tau^{(R)} - \eta)
\hat{\varepsilon}_1^{\frac{3}{2}} + O(\hat{\varepsilon}_1^2)\;.
\label{eq:PdmHB}
\end{equation}

\noindent
{\bf Step 6:}~~~{\em Obtain the full Poincar\'{e} map, $\mathcal{P}$ and compute $h_3$.}\\
The full Poincar\'{e} map is
$\mathcal{P} = \mathcal{P}_{\rm lhf} \circ \mathcal{P}_{\rm dm}$.
In scaled coordinates (\ref{eq:PlhfHB}) becomes
\begin{equation}
\hat{\varepsilon}_6 =
\eta_2 (\frac{\pi}{\omega^3} + O(|\mu,\eta_2|^1)) +
\left( 1 + \frac{\pi}{\omega} \eta_2 + q_1 \mu^2 + q_2 \mu \eta_2 +
q_3 \eta_2^2 + O(|\mu,\eta_2|^3) \right) \hat{\varepsilon}_4 +
O(\hat{\varepsilon}_4^2)\;.
\label{eq:PlhfHB2}
\end{equation}
Composing (\ref{eq:PdmHB}) and (\ref{eq:PlhfHB2}) produces
\begin{equation}
\hat{\varepsilon}_6 = \mathcal{P}(\hat{\varepsilon}_1;\mu,\eta_2)
= \Omega_0 + \Omega_1 \hat{\varepsilon}_1
+ \Omega_2 \hat{\varepsilon}_1^{\frac{3}{2}} + O(\hat{\varepsilon}_1^2)\;,
\label{eq:theMapP}
\end{equation}
where
\begin{equation}
\begin{split}
\Omega_0(\mu,\eta_2) & = \frac{\pi}{\omega^3} \eta_2 + O(2)\;, \\
\Omega_1(\mu,\eta_2) & = 1 + \frac{\pi}{\omega} \eta_2
+ q_1 \mu^2 + q_2 \mu \eta_2 + q_3 \eta_2^2 + O(3)\;, \\
\Omega_2(\mu,\eta_2) & = \frac{4 \sqrt{2}}{3} \tau_R + O(1)\;.
\end{split}
\nonumber
\end{equation}
To remove fractional powers we introduce $\chi = \sqrt{\hat{\varepsilon}}$.
The function
\begin{equation}
G_2(\chi;\mu,\eta_2) = \mathcal{P}(\chi^2;\mu,\eta_2) - \chi^2\;,
\nonumber
\end{equation}
is $C^{k-1}$, and by the implicit function theorem,
there exists a unique $C^{k-1}$ function, $\mathcal{F}$,
such that $G_2(\chi;\mu,\mathcal{F}(\chi;\mu)) = 0$.
Via a series expansion it is straightforward to obtain
\begin{equation}
\mathcal{F}(\chi;\mu) = \left( -\frac{\omega^3 q_1}{\pi} \mu^2 + O(\mu^3) \right) \chi^2 +
\left( -\frac{4 \sqrt{2} \omega^3 \tau_R}{3 \pi} + O(\mu) \right) \chi^3 + O(\chi^4)\;.
\label{eq:fpsHB}
\end{equation}
The fixed point (\ref{eq:fpsHB}), has an associated multiplier of one
when the $C^{k-2}$ function
\begin{eqnarray}
G_3(\chi;\mu) & = & \frac{\partial \mathcal{P}}{\partial \varepsilon}
(\chi^2;\mu,\mathcal{F}(\chi;\mu)) - 1 \nonumber \\
& = & (q_1 \mu^2 + O(\mu^3)) + (2\sqrt{2} \tau_R + O(\mu)) \chi + O(\chi^2)\;, \nonumber
\end{eqnarray}
is zero.
By the implicit function theorem,
there exists a unique $C^{k-2}$ function, $\hat{h}_3$, such that
$G_3(\hat{h}_3(\mu);\mu) = 0$.
Furthermore, using (\ref{eq:q1HB}),
\begin{equation}
\hat{h}_3(\mu)
= -\frac{\sqrt{2} \pi a_0}{\omega^5 \tau_R} \mu^2 + O(\mu^3)\;.
\nonumber
\end{equation}
Notice this fixed point is valid when
$\hat{\varepsilon}_1^{\frac{1}{2}} = \chi = \hat{h}_3(\mu) \ge 0$
which is true when $a_0$ and $\tau_R$ have opposite signs.
Finally let
\begin{eqnarray}
h_3(\mu) & = & h_2(\mu) + \mathcal{F}(\hat{h}_3(\mu);\mu) \nonumber \\
& = & h_2(\mu) - \frac{8 \pi^2 a_0^3}{3 \omega^{12} \tau_R^2} \mu^6 + o(\mu^6)\;. \nonumber
\end{eqnarray}
Then $h_3$ is the $C^{k-2}$ function (\ref{eq:h3HB})
and we have verified (iii) and (iv) of the theorem.
$\Box$

\section{Conclusions}
\label{sec:conclusions}
We have presented an unfolding of the codimension-two
simultaneous occurrence of a discontinuous bifurcation and an Andronov-Hopf bifurcation
in a general, piecewise-smooth, continuous system.
We have found a locus of Hopf bifurcations that emanates from the
codimension-two point, (\ref{eq:h1HB}).
Tangent to this is a locus of grazing bifurcations of the Hopf cycle
with the switching manifold, (\ref{eq:h2HB}).
Heuristically, the curves are tangent because,
with respect to a linear change in parameter values,
the distance between the equilibrium solution and the switching manifold
increases linearly, whereas the amplitude of the Hopf cycle grows
as the square root of the magnitude of the parameter change.

A periodic orbit is created at the discontinuous bifurcation
on one side of the codimension-two point.
When the stability of this orbit opposes that of the Hopf cycle,
the two orbits collide and annihilate in a saddle-node bifurcation
on a curve that deviates only to order six from the grazing bifurcation (\ref{eq:h3HB}).
The mechanism behind this sixth order scaling
can be explained with a simple calculation.
Omitting higher order terms,
the Poincar\'{e} map, $\mathcal{P}$, (\ref{eq:theMapP}),
is essentially
\begin{equation}
\varepsilon' = \eta_2 + \Xi(\mu) \varepsilon +
\gamma \varepsilon^{\frac{3}{2}}\;,
\label{eq:simpleMap}
\end{equation}
where $\Xi(\mu)$ is the Floquet multiplier of the Hopf cycle
and $\gamma$ is a constant.
The grazing bifurcation occurs when $\eta_2 = 0$.
It is easily determined, as in \cite{DiBu08},
a saddle-node bifurcation of the fixed point of (\ref{eq:simpleMap}) occurs when
\begin{equation}
\eta_2^{\rm SN} = \frac{4}{27 \gamma^2}(1 - \Xi)^3\;.
\end{equation}
The Floquet multiplier is unity at the Hopf bifurcation ($\eta_1 = 0$)
and by assumption varies linearly with respect to $\eta_1$,
i.e., $\Xi \approx 1 + \lambda \eta_1$, ($\lambda \ne 0$).
The grazing bifurcation occurs when $\eta_1 = O(\mu^2)$,
thus $\Xi \approx 1 + \hat{\lambda} \mu^2$, ($\hat{\lambda} \ne 0$).
Hence
\begin{equation}
\eta_2^{\rm SN} \approx -\frac{4 \hat{\lambda}^3}{27 \gamma^2} \mu^6\;.
\end{equation}

For simplicity, throughout this paper we assumed the switching manifold
was infinitely differentiable. 
Our analysis is unchanged if the switching manifold is only $C^3$.
However, if the switching manifold were $C^2$ and not $C^3$,
we would be unable to determine the same expression for the map
$\mathcal{P}_3$, (\ref{eq:P3HB}).

Recently it has been found that an eight-dimensional model of
yeast growth \cite{SiKo08} exhibits codimension-two
discontinuous bifurcations such as the scenario described here.
Other observed codimension-two situations that remain to be rigorously unfolded
include the simultaneous
occurrence of a saddle-node and discontinuous bifurcation,
and a discontinuous Hopf bifurcation \cite{FrPo97,SiMe07}
of indeterminable criticality.
We hope to report on these in a future paper, see also \cite{Si09}.

We would also like to extend the results of this paper to
higher dimensional systems, like the yeast model.
It seems plausible that bifurcation sets for higher-dimensional systems
will exhibit scalings of the same orders, but we do not, as of yet,
have a formal justification of this.

\end{document}